\documentclass[12pt]{article}
\usepackage{amsmath}
\usepackage{amssymb}
\usepackage{amscd}

\newtheorem{Df}{Definition}[section]
\newtheorem{Te}[Df]{Theorem}
\newtheorem{Po}[Df]{Proposition}
\newtheorem{Cr}[Df]{Corollary}
\newtheorem{Lm}[Df]{Lemma}
\newtheorem{Ca}[Df]{Claim}
\newtheorem{Cn}[Df]{Conjecture}

\newcommand{\Bdf}{\begin{Df}}
\newcommand{\Edf}{\end{Df}}
\newcommand{\Bte}{\begin{Te}}
\newcommand{\Ete}{\end{Te}}
\newcommand{\Bpo}{\begin{Po}}
\newcommand{\Epo}{\end{Po}}
\newcommand{\Bcr}{\begin{Cr}}
\newcommand{\Ecr}{\end{Cr}}
\newcommand{\Blm}{\begin{Lm}}
\newcommand{\Elm}{\end{Lm}}
\newcommand{\Bca}{\begin{Ca}}
\newcommand{\Eca}{\end{Ca}}
\newcommand{\Bcn}{\begin{Cn}}
\newcommand{\Ecn}{\end{Cn}}
\newcommand{\Bdm}{{\it Proof.}\ }
\newcommand{\Edm}{\rule{2mm}{2mm}}
\newcommand{\Rm}{{\it Remark \arabic{section}.\arabic{Df} \ }}
\newcommand{\Ea}{{\it Example \arabic{section}.\arabic{Df} \ }}

\DeclareMathOperator{\Hom}{\mathrm{Hom}}
\newcommand{\G}{\Gamma}

\newcommand{\eps}{\epsilon}
\newcommand{\si}{\sigma}
\newcommand{\Si}{{\mathbb{S}}}

\newcommand{\beq}{\begin{equation}\label}
\newcommand{\eeq}{\end{equation}}
\newcommand{\qed}{\hfill$\Box$\medskip}
\newcommand{\op}{\operatorname}
\newcommand{\Id}{\op{Id}}
\newcommand{\map}{\longrightarrow}
\newcommand{\sset}{\subset}
\renewcommand{\o}{\otimes}

\begin{document}

\title{\bf{Symplectic reflection algebras and non-homogeneous
    $N$-Koszul property}}
\author{\bf{Roland Berger} \textrm{and} \bf{Victor Ginzburg}\\ \\ 
\emph{E-mail}: Roland.Berger@univ-st-etienne.fr \\ ginzburg@math.uchicago.edu}
\date{}
\maketitle
\begin{abstract}
From symplectic reflection algebras~\cite{eg:reflec}, some algebras
are naturally introduced. We show that these algebras are non-homogeneous
$N$-Koszul algebras. The Koszul property was generalized to homogeneous algebras of degree $N>2$
in~\cite{rb:nonquad}. In the present paper, the extension of the Koszul property
to non-homogeneous algebras is realized through a PBW theorem. This
PBW theorem is the generalization to the $N$-case of a quadratic result
obtained by Braverman and Gaitsgory~\cite{bg:pbw} (see also Polishchuk and
Positselski~\cite{pp:quad}). Symplectic reflection algebras need to
work over special non-commutative semi-simple rings. Actually,
replacing the ground field by any von Neumann regular ring is more
general and is well-adapted to the Koszul property and the PBW theorem. 
\end{abstract}

\section{Introduction}
 Let $V$ be a finite dimensional complex vector space which is endowed with a symplectic 2-form $\omega$. Let $\G$ be a finite subgroup of Sp$(V)$ and $(TV)\# \G$ the smash product of the tensor algebra $TV$ of $V$ with the group algebra $\mathbb{C}\G$ of $\G$. For any $g \in \G$, introduce the subspaces $M_g:=\op{Im}(\Id- g)$ and $L_g:=\op{Ker}(\Id- g)$, so that one has $V=M_g \oplus L_g$ and 
\begin{equation} \label{deco}
\Lambda ^{2}(V)=(\Lambda ^2(M_g)) \oplus (M_g\otimes L_g) \oplus (\Lambda ^2(L_g)).
\end{equation}

Recall that the integer $a(g):=\dim M_g$ is even, and that $g$ is called a \emph{symplectic reflection} if this dimension is 2. Define the $\mathbb{C}$-linear map $\psi_g : \Lambda ^{2}(V) \rightarrow \mathbb{C}$ in order to coincide with $\omega$ on $\Lambda ^{a(g)}(M_g)\otimes \Lambda ^{2-a(g)}(L_g)$ and to vanish on the other components of (\ref{deco}). Clearly, $\psi_g=0$ if $g$ is neither $Id$ nor a symplectic reflection. 

Writing $\psi=\sum_{g\in\G}\psi_g\cdot g$, we define a skew-symmetric $\mathbb{C}$-bilinear pairing $\psi : V\times V \rightarrow \mathbb{C} \G$ which is $\G$-equivariant. Then the \emph{symplectic reflection algebra}~\cite{eg:reflec} is the $\mathbb{C} \G$-algebra H$_{\psi}$ defined by
$$H_{\psi} = (TV\# \Gamma) / I(x\otimes y-y \otimes x-\psi (x,y)\, ; \ x,\, y \in V),$$
Actually, for any map $m:\G \rightarrow \mathbb{C}$ which is constant on any conjugation class, $m\cdot \psi$ is $\G$-equivariant, and H$_{m\cdot \psi}$ is also called a symplectic reflection algebra. 

The algebra H$_{m\cdot \psi}$ is naturally filtered, and there is a natural graded algebra morphism $H_{0\cdot \psi}= (SV)\# \G \rightarrow gr(H_{m\cdot \psi})$. The first fundamental feature of H$_{m\cdot \psi}$ is that this morphism is an \emph{isomorphism} : it is the so-called PBW property for H$_{m\cdot \psi}$. The aim of this paper is to show that such a property makes sense in a more general context than the quadratic one, and that the generalized PBW property holds for a new class of H$_{m\cdot \psi}$. Let us describe briefly this new class.

Now $p$ is an integer with $2\leq p \leq \dim V$, $\G$ is a finite subgroup of GL$(V)$, and $\phi : \Lambda^pV\to k$ is a $\G$-invariant linear map (playing the role of $\omega$). We have an analogous decomposition (\ref{deco}) for $\Lambda ^{p}(V)$, and analogous definitions for the $\psi_g$'s, $\psi$, $m\cdot \psi$, and H$_{m\cdot \psi}$. For the relations of H$_{m\cdot \psi}$, we replace the 2-tensor $x\otimes y-y \otimes x$ by any totally skew-symmetric tensor of $p$ variables, while $m\cdot \psi$ is applying on these variables. We shall prove the following result (Corollary 4.5 below) which states a PBW property for the new class.
\Bte \label{intro1}
Under the previous notations and assumptions, the natural graded algebra morphism $H_{0\cdot \psi} \rightarrow gr(H_{m\cdot \psi})$ is an isomorphism.
\Ete

The PBW property in the case of any filtered algebra with (inhomogeneous) quadratic relations is well understood in the setting due to Braverman and Gaitsgory~\cite{bg:pbw} (see also Polishchuk and Positselski~\cite{pp:quad}). This setting brings out the fundamental role of the Koszul property (in Priddy's sense) of the homogeneous quadratic algebra which is obtained by forgetting the non-quadratic part. 

The first author has extended the Koszul property to any algebra with $N$-homogeneous relations ($N$ is fixed $\geq 2$)~\cite{rb:nonquad}. So an analog of the PBW theorem for any filtered algebra with $N$-inhomogeneous relations is naturally the first step for our proof of Theorem \ref{intro1}. We shall prove such a PBW theorem under very large assumptions, including the fact that the ground field $\mathbb{C}$ has to be replaced by the group algebra $\mathbb{C} \G$ which is a non-commutative ring! Fortunately, this ring is semi-simple, and our formalism for Koszul and PBW properties works out even for more general rings, the von Neumann regular rings. Notice that the extension of the Koszul property to any semi-simple ground ring was already performed in~\cite{bgso:kdp} in the quadratic case. Our PBW theorem is the following (Theorem 3.4 below).
\Bte \label{intro2}
\emph{(PBW theorem in the $N$-case)} Assume that $k$ is a von Neumann regular ring, $V$ is a
$k$-$k$-bimodule, $N$ is an integer $\geq 2$, and $P$ is a sub-$k$-$k$-bimodule of $F^N$, where
$F^n=\oplus_{0\leq i \leq n} V^{\otimes i}$ for any $n \geq 0$. Set $U=T(V)/I(P)$ and
$A=T(V)/I(R)$, where $R=\pi(P)$ and $\pi$ is the projection of $F^N$ onto $V^{\otimes N}$ modulo
$F^{N-1}$. 

Assume that the graded left $k$-module Tor$_3^A(k_A,_Ak)$ is concentrated in degree $N+1$ (this property is a part of the Koszul property of $A$).
Then the conditions
$$P\cap F^{N-1} =0,$$
$$(PV+VP) \cap F^{N} \subseteq P,$$
imply that the PBW property holds, i.e., the natural algebra morphism
$A\rightarrow gr(U)$ is an isomorphism.
\Ete

This theorem suggests to extend the terminology ``$A$ is Koszul'' to ``$U$ is Koszul'' (see Definition 3.9 below). It is natural since J. L. Koszul introduced his resolution for the polynomial algebra  in the \emph{filtered} context of the enveloping algebra of a Lie algebra, and used the classical PBW property as a trick to carry over the exactness of his resolution to the standard complex~\cite{weib:homo}.

\setcounter{equation}{0}

\section{$N$-Koszul algebras over von Neumann regular rings}
A ring $k$ is said to be \emph{von Neumann regular} if for every $x\in k$, there exists
$y\in k$ such that $xyx=x$. In this text, von Neumann regular rings will be used only
through the following characterization (in which left can be replaced by right)~\cite{good:von}.
\Bpo
A ring $k$ is von Neumann regular if and only if all left $k$-modules are flat.
\Epo
Von Neumann regular rings are exactly the rings having weak dimension 0~\cite{weib:homo}. So a
semi-simple ring (i.e., a ring having left or right global dimension 0) is von Neumann regular. An
infinite product
$k$ of fields is von Neumann regular, but is not semi-simple. The same holds for $k=End(V)$, where
$V$ is an infinite-dimensional vector space.

Throughout this section, $k$ is a von Neumann regular ring, and $V$ denotes a graded
$k$-$k$-bimodule which is \emph{concentrated in degree 1}. The tensor power (over $k$)
$V^{\otimes n}$ for $n= 0,1,\ldots$ is a graded $k$-$k$-bimodule which is concentrated in degree $n$.
The direct sum $T(V)=\bigoplus_{n\geq 0}V^{\otimes n}$ is naturally a
$\mathbb{N}$-graded ring, and a $k$-$k$-bimodule whose left or right actions coincide with the
products in the ring by elements of
$V^{\otimes 0}=k$. We sum up these properties by saying that $T(V)$ is a
connected $\mathbb{N}$-graded $k$-$k$-algebra.

If $a$ and $b$ are in $T(V)$, their product in $T(V)$ is denoted by $ab$. For any
sub-$k$-$k$-bimodules $E$ and $F$ of $T(V)$, $EF$ denotes the sub-$k$-$k$-bimodule formed by finite
sums of products $ab$, $a\in E$, $b\in F$. The two-sided ideal $I(E)$ of $T(V)$ generated by $E$ is
such that
$$I(E)=T(V)E\,T(V)= \sum_{i,j\geq 0} V^{\otimes i}E\,V^{\otimes j}.$$
On the other hand, if $E\subseteq V^{\otimes i}$ and $F\subseteq V^{\otimes
j}$, the canonical homomorphism $E\otimes_k F\rightarrow V^{\otimes (i+j)}$ is injective since the
right $k$-module $E$ and the left $k$-module $F$ are \emph{flat}. So the $k$-$k$-bimodules
$E\otimes_k F$ and $EF$ will be identified. The following result is known (\cite{bou:algco},
chap.I, ¤2, n$^{\circ}$ 6) and is again a consequence of flatness.
\Blm \label{formulas}
If $E$ and $E'$ are sub-$k$-$k$-bimodules of $V^{\otimes i}$, and if $F$ and $F'$ are
sub-$k$-$k$-bimodules of $V^{\otimes j}$, the following formulas hold :

(i) $(EF)\cap (EF')=E(F\cap F')$,

(ii) $(EF)\cap (E'F)=(E\cap E')F$,

(iii) $(E'F)\cap (EF')=E'F'$ if moreover $E'\subseteq E$ and $F'\subseteq F$.
\Elm

Now fix an integer $N\geq 2$ and a sub-$k$-$k$-bimodule $R$ of $V^{\otimes N}$. Then the
two-sided ideal $I(R)$ is graded by the $k$-$k$-bimodules
$$I(R)_n=\sum_{i+N+j=n} V^{\otimes i}R\, V^{\otimes j}, \ n\geq 0.$$
One has $I(R)_n=0$ if $0\leq n \leq N-1$. In this section, we are interested in the connected
$\mathbb{N}$-graded $k$-$k$-algebra $A=T(V)/I(R)$. The gradation of $A$ is formed by the
$k$-$k$-bimodules $A_n=V^{\otimes n}/I(R)_n$. One has $A_n=V^{\otimes n}$ if $0\leq n \leq N-1$. 
The natural projection $\epsilon : A \rightarrow A_0=k$ makes $k$ as being an $A$-$A$-bimodule,
denoted by $_Ak_A$. We also use $_Ak$ and $k_A$ for the left and right associated $A$-modules.

The categories of left $A$-modules, right $A$-modules, $A$-$A$-bimodules are respectively denoted by
$A$-Mod, Mod-$A$, $A$-Mod-$A$. When the objects are graded and the arrows are homogeneous of degree
0, the categories are respectively denoted by $A$-grMod, grMod-$A$, $A$-grMod-$A$. For example, the
object $_Ak$ is graded (as concentrated in degree 0), so that $\epsilon$ is an arrow of
$A$-grMod. Any object $M$ of $A$-grMod has a free resolution in $A$-grMod, i.e., a resolution by
graded-free left $A$-modules. Accordingly, for any natural number $n$, Tor$_n^A(k_A,M)$ is a
graded left $k$-module. Koszul property of $A$ will be defined from the objects
Tor$_n^A(k_A,_Ak)$, $n\in \mathbb{N}$. Let us begin to compute these objects for $n=0,1,2$.

Actually, for computing Tor$_n^A(k_A,_Ak)$, it is enough to have a \emph{flat} resolution of $_Ak$
in $A$-grMod, i.e., a resolution in $A$-grMod formed by flat left $A$-modules (an object of
$A$-grMod which is flat in $A$-Mod is flat in $A$-grMod). On the other hand, it is well-known
(and easy to prove) that, for any homomorphism from a ring $k$ to a ring $A$ and any flat left
$k$-module $E$, the left $A$-module $A\otimes_k E$ is flat (\cite{bou:algco},
chap.I, ¤2, n$^{\circ}$ 7). Thus, since $k$ is von Neumann regular, $A\otimes_k E$ is
a \emph{flat} left $A$-module (not free in general!) for any left $k$-module $E$. So we search
resolutions in $A$-grMod whose objects are of the type $A\otimes_kE$, $E$ object of
$k$-grMod.
\Blm \label{ident}
Fix the natural numbers $n$ and $m$, and let $E$ be a sub-$k$-$k$-bimodule of $V^{\otimes m}$. The
$k$-$k$-bimodules $A_n\otimes_k E$ and $V^{\otimes n}E / I(R)_n E$ are naturally isomorphic.
\Elm
\Bdm
As $E$ is flat in $k$-Mod, one has the natural exact sequence of $k$-$k$-grMod
$$0\rightarrow I(R)_n\otimes_k E \rightarrow V^{\otimes n}\otimes_k E \rightarrow A_n\otimes_k E
\rightarrow 0.$$
But $V^{\otimes n}\otimes_k E$ is identified to the sub-$k$-$k$-bimodule $V^{\otimes n}E$ of
$T(V)$, and the image of the injective map is identified to the sub-$k$-$k$-bimodule $I(R)_n E$ of
$T(V)$.
\Edm
\\ 

The arrow $\epsilon : A\rightarrow \, _Ak$ of $A$-grMod has $A_{\geq 1}=\bigoplus_{n\geq 1} A_n$ as
kernel. The inclusion $V \rightarrow A$ defines an injective (by flatness!) natural arrow
$A\otimes_kV \rightarrow A\otimes_kA$, which is composed with the multiplication $\mu : A\otimes_kA
\rightarrow A$ to define an arrow $\delta_1:A\otimes_k V \rightarrow A$ of $A$-grMod. In degree 0,
$\delta_1$ vanishes, and in degree $n\geq 1$, it is identified by Lemma \ref{ident} to the
canonical map
$$\frac{V^{\otimes n}}{I(R)_{n-1}V} \rightarrow \frac{V^{\otimes n}}{I(R)_n},$$
which is surjective. So im$(\delta_1)=\ker(\epsilon)$.

The natural injections $R\rightarrow V^{\otimes N} = V^{\otimes(N-1)} \otimes _k V \rightarrow
A\otimes_k V$ define a natural injection $A\otimes_k R \rightarrow A\otimes_k(A\otimes_k V)\cong
(A\otimes_kA)\otimes_k V$, which is composed with $\mu\otimes_k 1_V$ to define an arrow $\delta_2 :
A\otimes_k R \rightarrow A\otimes_k V$ of $A$-grMod. In degrees $<N$, $\delta_2$ vanishes, and in
degree $n\geq N$, it is identified to the canonical map
\begin{equation} \label{degree2}
\frac{V^{\otimes (n-N)}R}{I(R)_{n-N}R} \rightarrow \frac{V^{\otimes n}}{I(R)_{n-1}V}\ ,
\end{equation}
so that one has
$$(\mathrm{im}(\delta_2))_n=\frac{V^{\otimes (n-N)}R+ I(R)_{n-1}V}{I(R)_{n-1}V}\ .$$
On the other hand, for $n\geq 1$, one has
$$(\ker (\delta_1))_n= \frac{I(R)_n}{I(R)_{n-1}V}\ .$$
Comparing the two equalities, we get im$(\delta_2)=\ker (\delta_1)$. Clearly, $\ker(\delta_2)$
vanishes in degrees $\leq N$. Choose an object $E_3$ of $k$-grMod living in degrees $\geq N+1$,
and a surjective arrow $E_3\rightarrow
\ker(\delta_2)$ in $k$-grMod. Extending the latter by $A$-linearity, one gets a surjective arrow
$A\otimes_k E_3\rightarrow \ker(\delta_2)$, which is composed with the inclusion
$\ker(\delta_2)\rightarrow A\otimes_k R$ to define an arrow $\delta_3 : A\otimes_k E_3 \rightarrow
A\otimes_k R$ of $A$-grMod. Finally
\begin{equation} \label{beginning}
A\otimes_k E_3 \stackrel{\delta_{3}}{\longrightarrow} A\otimes_k R
\stackrel{\delta_{2}}{\longrightarrow} A\otimes_k V \stackrel{\delta_{1}}{\longrightarrow} A
\longrightarrow 0
\end{equation}
is the beginning of a flat resolution of $_Ak$ (via $\epsilon :A \rightarrow k$) in $A$-grMod. 

Applying the functor $k_A \otimes_A -$ to (\ref{beginning}), the sequence 
$$E_3 \longrightarrow R \longrightarrow  V \longrightarrow k \longrightarrow 0$$ 
of $k$-grMod is obtained, in which all arrows are vanishing (remember that the arrows in $k$-grMod
are of degree 0). Therefore, one has the isomorphisms in $k$-grMod : Tor$_{0}^{A}(k_A,_Ak)\cong k$, 
Tor$_{1}^{A}(k_A,_Ak)\cong V$, Tor$_{2}^{A}(k_A,_Ak)\cong R$. Moreover Tor$_3^A(k_A,_Ak)$ lives in
degrees $\geq N+1$. 

The three first Tor's are concentrated in only one degree (0, 1, $N$ respectively).
Roughly speaking, the Koszul property means that each Tor$_n^A(k_A,_Ak)$ is concentrated in the
lowest possible degree when $n=0,1,2,3,\ldots$ We have to examine now when Tor$_3^A(k_A,_Ak)$ is
concentrated in degree $N+1$.
\Bpo
The graded left $k$-module Tor$_3^A(k_A,_Ak)$ is concentrated in degree $N+1$ if and only if the
graded left $A$-module $\ker (\delta_2)$ is generated in degree $N+1$.
\Epo
\Bdm
By dimension shifting applied to flat resolutions (\cite{weib:homo}, p.47), one has
$$\mathrm{Tor}_3^A(k_A,_Ak)= \ker \left(k_A\otimes_A \ker(\delta_2) \stackrel{k_A\otimes_A
i}{\longrightarrow} k_A\otimes_A(A\otimes_kR)\cong R \right),$$
where $i: \ker(\delta_2) \rightarrow A\otimes_k R$ is the inclusion. Since $\ker (\delta_2)$ lives
in degrees $\geq N+1$, $k_A\otimes_A \ker(\delta_2)$ lives in degrees $\geq N+1$. But $R$ is
concentrated in degree $N$ and $k_A\otimes_A i$ preserves the degrees. Thus
Tor$_3^A(k_A,_Ak)=k_A\otimes_A\ker(\delta_2)$. It suffices to prove the following. 
\Blm
Fix $n \in \mathbb{Z}$. Let $M=\bigoplus_{i\geq n} M_i$ be an object of $A$-grMod, living in degrees
$\geq n$. Then $M$ is generated in degree $n$ if and only if the graded left $k$-module
$k_A\otimes_AM$ is concentrated in degree $n$. 
\Elm
\Bdm
The graded left $k$-module $k_A\otimes_AM$ is canonically isomorphic to $M/A_{\geq 1}.M$. Assuming
$M=A.M_n$, one sees that $M_i=0$ if $i<n$, and for $i>n$, $M_i=A_{i-n}.M_n$ is contained in $(A_{\geq
1}.M)_i$, hence $M/A_{\geq 1}.M$ vanishes in degrees $\neq n$. Conversely, assuming that
$M_i=(A_{\geq 1}.M)_i$ for $i>n$, one has $M_{n+1}=A_1.M_n$, and inductively $M_i=A_{i-n}.M_n$ for
any $i>n$, so $M=A.M_n$. 
\Edm
\\

Identifying $\delta_2$ in degree $n$ with (\ref{degree2}) provides 
\begin{equation} \label{ker2}
(\ker(\delta_2))_n= \frac{(V^{\otimes (n-N)}R)\cap (I(R)_{n-1}V)}{I(R)_{n-N}R}\ ,\, n\geq N+1.
\end{equation}
It follows
$$(\ker(\delta_2))_n=(V^{\otimes (n-N)}R)\cap (R V^{\otimes (n-N)}+\cdots + V^{\otimes
(n-N-1)}RV)\ ,\, N+1\leq n \leq 2N-1,$$
and in particular $(\ker(\delta_2))_{N+1}=W_{N+1}$, where the following notation is used
$$W_n=\bigcap_{i+N+j=n}V^{\otimes i}RV^{\otimes j}\ ,\, n\geq N.$$
In our choice of the graded left $k$-module $E_3$, we can assume that $(E_3)_{N+1}=W_{N+1}$, so
that $\delta_3$ is injective in degree $N+1$. Thus $\ker(\delta_3)$ lives in degrees $\geq N+2$, and
$$(\mathrm{Tor}_3^A(k_A,_Ak))_{N+1}=W_{N+1}.$$

Coming back to our initial question, we have to know when $\ker(\delta_2)$ is generated in degree
$N+1$. It is equivalent to saying that for every $n\geq N+2$, one has
\begin{equation} \label{cond}
(\ker(\delta_2))_n= A_{n-N-1}.W_{N+1}.
\end{equation}
So a partial condition is obtained for $N+2\leq n \leq 2N-1$ :
\begin{equation} \label{ec}
(V^{\otimes (n-N)}R)\cap (R V^{\otimes (n-N)}+\cdots + V^{\otimes (n-N-1)}RV)=V^{\otimes
(n-N-1)}W_{N+1}.
\end{equation}
The set of equalities (\ref{ec}) when $N+2\leq n \leq 2N-1$ is called the \emph{extra condition},
and is denoted by (ec). That condition does not occur if $N=2$ (hence its name), and it reveals the
jump of degrees between generators and relations. It is remarkable that (ec) will be sufficient to
take into account the other jumps of degrees which will appear in the inductive process.

Examine now (\ref{cond}) for $n\geq 2N$. Since $W_{N+1}$ is a left sub-$k$-module of $V\otimes_k R$,
$A.W_{N+1}$ is the left sub-$A$-module of $A.(V\otimes_k R)\subseteq A\otimes_k R$ which is the image
of the composite
$$A\otimes_k W_{N+1} \stackrel{1_A \otimes i}{\longrightarrow} A\otimes_k V \otimes_k R
\stackrel{\delta_{1}\otimes 1_R}{\longrightarrow} A\otimes_k R,$$
where $i: W_{N+1}\rightarrow V\otimes_k R$ is the inclusion. Using Lemma \ref{ident}, the maps of
this composite become in each degree $n\geq N+1$, the canonical following ones
$$\frac{V^{\otimes (n-N-1)}W_{N+1}}{I(R)_{n-N-1}W_{N+1}} \rightarrow \frac{V^{\otimes
(n-N)}R}{I(R)_{n-N-1}VR} \rightarrow \frac{V^{\otimes (n-N)}R}{I(R)_{n-N}R}\ .$$
Therefore 
$$A_{n-N-1}.W_{N+1}= \frac{V^{\otimes (n-N-1)}W_{N+1}+I(R)_{n-N}R}{I(R)_{n-N}R}\ ,$$
and joining that with (\ref{ker2}), we get the answer to our initial question.
\Bte \label{tor3}
The left $k$-module Tor$_3^A(k_A,_Ak)$ is concentrated in degree $N+1$ if and only if (ec)
and the following relations
\begin{equation} \label{other}
(V^{\otimes (n-N)}R)\cap (I(R)_{n-1}V)=V^{\otimes (n-N-1)}W_{N+1}+I(R)_{n-N}R\ ,\, n\geq 2N,
\end{equation}
are satisfied.
\Ete
A triple $(E,F,G)$ of sub-$k$-$k$-bimodules of $T(V)$ is said to be \emph{distributive} if $E\cap
(F+G)=E\cap F +E\cap G$, and the latter equality is called a \emph{distributivity relation}. Show
that, if (ec) holds, (\ref{other}) is a distributivity relation for each $n$. In fact, setting
$$E= V^{\otimes (n-N)}R,\ F= I(R)_{n-N}V^{\otimes N},\ G=V^{\otimes (n-2N+1)}I(R)_{2N-2}V,$$
the left-hand side of (\ref{other}) is $E\cap(F+G)$. Using Lemma \ref{formulas}, (iii), one has
$$E\cap F=I(R)_{n-N}R.$$
The part (i) of this lemma provides
$$E\cap G= V^{\otimes (n-2N+1)}[(V^{\otimes (N-1)}R)\cap (R V^{\otimes (N-1)}+\cdots + V^{\otimes
(N-2)}RV)],$$
and (ec) for $n=2N-1$ implies
$$E\cap G= V^{\otimes (n-2N+1)}[V^{\otimes (N-2)}W_{N+1}]= V^{\otimes (n-N-1)}W_{N+1},$$
so the right-hand side of (\ref{other}) is $E\cap F +E\cap G$.

Theorem \ref{tor3} will be enough to state the PBW theorem, but we want to give further indications
about the definition of the Koszul property. The details are left to the reader. Assume now that
Tor$_3^A(k_A,_Ak)$ is concentrated in degree $N+1$. In particular, (ec) holds. The 
inclusion $W_{N+1}\rightarrow A\otimes_k R$ defines a natural injection $A\otimes_k W_{N+1}
\rightarrow A\otimes_k(A\otimes_k R)\cong (A\otimes_kA)\otimes_k R$, which is composed with
$\mu\otimes_k 1_R$ to define a more suitable arrow $\delta_3 : A\otimes_k W_{N+1} \rightarrow
A\otimes_k R$ of $A$-grMod. The beginning of the flat resolution is now
\begin{equation} \label{beginning2}
A\otimes_k W_{N+1} \stackrel{\delta_{3}}{\longrightarrow} A\otimes_k R
\stackrel{\delta_{2}}{\longrightarrow} A\otimes_k V \stackrel{\delta_{1}}{\longrightarrow} A
\longrightarrow 0.
\end{equation}
In degrees $<N+1$, $\delta_3$ vanishes, and for $n\geq N+1$, one has
$$(\ker(\delta_3))_n= \frac{(V^{\otimes (n-N-1)}W_{N+1})\cap (I(R)_{n-N}R)}{I(R)_{n-N-1}W_{N+1}}\
.$$ 
Clearly, $(\ker(\delta_3))_n=0$ if $N+1\leq n \leq 2N-1$, and 
$$(\ker(\delta_3))_{2N}=(V^{\otimes (N-1)}W_{N+1})\cap (RR).$$
Lemma \ref{formulas} (iii) implies that $RR=(RV^{\otimes N})\cap (V^{\otimes N}R)$. On the
other hand, it is easy to check that (ec) implies that 
$$(V^{\otimes (N-1)}R)\cap (RV^{\otimes (N-1)})=W_{2N-1}.$$
So, using Lemma \ref{formulas} (ii), we get
$(\ker(\delta_3))_{2N}=W_{2N}$. In other words, Tor$_4^A(k_A,_Ak)$ lives in degrees $\geq 2N$ and 
$$(\mathrm{Tor}_4^A(k_A,_Ak))_{2N}=W_{2N}.$$
The next question would be to know when Tor$_4^A(k_A,_Ak)$ is concentrated in degree $2N$.

More generally, in our inductive process to define the Koszul property, the successive degrees of the
Tor's are given by the values of the so-called \emph{jump map} $\zeta : \mathbb{N}\rightarrow
\mathbb{N}$, where
$$\zeta(2q)=qN, \ \zeta(2q+1)= qN+1,\ q\ \mathrm{integer} \geq 0.$$ 
Let $n$ be $\geq 3$. Assume that for any $3\leq i\leq n$, Tor$_i^A(k_A,_Ak)$ is
concentrated in degree $\zeta (i)$. In particular, (ec) holds. Then Tor$_{n+1}^A(k_A,_Ak)$ lives in
degree $\geq \zeta (n+1)$, and Tor$_{n+1}^A(k_A,_Ak)$ is concentrated in degree $\zeta (n+1)$ if
and only if a certain sequence of distributivity relations hold. This sequence of distributivity
relations is written down in Theorem 2.11 (iii) of~\cite{rb:nonquad}, where the case of a field $k$
was treated (but the statement holds for any von Neumann regular ring thanks to Lemma
\ref{formulas}). Anyway, the following definition becomes natural.
\Bdf
The connected $\mathbb{N}$-graded $k$-$k$-algebra $A$ is said to be Koszul if for any $n
\geq 3$, the graded left $k$-module Tor$_n^A(k_A,_Ak)$ is concentrated in degree $\zeta (n)$.
\Edf
When $A$ is Koszul, the inductive process constructs a flat resolution in $A$-grMod of $_Ak$,
which is called the \emph{Koszul resolution}. Actually, independently of the Koszul property of $A$, a
\emph{Koszul complex} $\mathcal{C}$ in $A$-grMod is easily constructed : the objects are the
$A\otimes_k
 W_{\zeta (n)}$, $n\geq 0$, and the arrow $A\otimes_k W_{\zeta (n+1)} \rightarrow A\otimes_k
W_{\zeta (n)}$ is the extension by left $A$-linearity of the inclusion $W_{\zeta (n+1)} 
\rightarrow A\otimes_k W_{\zeta (n)}$. If $A$ is Koszul, $\mathcal{C}$ is the Koszul resolution.
Conversely, if the complex $\mathcal{C}$ is exact in any (homological) degree $>0$, then
$\mathcal{C}$ is a resolution of $_Ak$ via $\epsilon : A
\rightarrow k$, and $A$ is Koszul since all arrows of the complex
$k_A\otimes_A \mathcal{C}$ are vanishing.

We shall need a change of rings result. Let $k'$ be a ring and $k\rightarrow k'$ a ring morphism.
Assume that the $k$-$k$-bimodule $V$ is a $k'$-$k$-bimodule whose left action extends left action by
$k$ and right action is the same. Assume also that $R$ is a left sub-$k'$-module of $V^{\otimes
N}$, so that $A$ is a $k'$-$k$-bimodule. Set $V'= V\otimes_k k'$. One has the sequence of
$k'$-$k'$-bimodule morphisms
$$V'\otimes_{k'} V' \cong (V\otimes_k k')\otimes_{k'} V' \cong
V\otimes_k (k'\otimes_{k'} V') \cong
V\otimes_k(V\otimes_k k') \cong (V\otimes_k V)\otimes_k k'.$$
More generally, $V^{\otimes_k n}\otimes_k k' \cong V^{'\otimes_{k'}
  n}$ for any $n$. As $k$ is von Neumann regular, $V^{\otimes_k
  n}\otimes_k k'$ is considered as included in $T_k(V)\otimes_k k'$, and
$T_k(V)\otimes_k k'$ is the direct sum of the $V^{\otimes_k
  n}\otimes_k k'$, $n\geq 0$. So we have a
natural isomorphism of $\mathbb{N}$-graded $k'$-$k'$-algebras
$$T_k(V)\otimes_k k' \cong T_{k'} (V')$$
which sends $I(R)\otimes_k k'$ (considered as included in
$T_k(V)\otimes_k k'$) onto $I(R')$, where $R'$ is the image of $R\otimes_k k'$.
The connected $\mathbb{N}$-graded $k'$-$k'$-algebra
$$A'= \frac{T_{k'} (V')}{I(R')}$$
is then canonically isomorphic to $A\otimes_k k'$ since 
$$\frac{T_k(V)\otimes_k k'}{I(R)\otimes_k k'} \cong \frac{T_k(V)}{I(R)}\otimes_k k'.$$ 
\Bpo \label{change}
With above notations and assumptions, assume that $k'$ is von Neumann regular and that the
connected $\mathbb{N}$-graded $k$-$k$-algebra $A$ with $N$-homogeneous relations is Koszul. Then
the connected $\mathbb{N}$-graded $k'$-$k'$-algebra $A'$ has $N$-homogeneous relations and is Koszul.
\Epo 
\Bdm
The Koszul complex $\mathcal{C}$ of $A$ is exact in degree $>0$. Then the
complex $\mathcal{C}\otimes_k k'$ is exact in degree $>0$, and its objects are naturally the
$A\otimes_k W'_{\zeta (n)}$, $n\geq 0$, where $W'_{\zeta (n)}\cong W_{\zeta (n)} \otimes_k k'$. But
$$A\otimes_k W'_{\zeta (n)} \cong A \otimes_k (k'\otimes_{k'}W'_{\zeta (n)}) \cong (A\otimes_k k')
\otimes_{k'}W'_{\zeta (n)} \cong A'\otimes_{k'}W'_{\zeta (n)},$$
so that $\mathcal{C}\otimes_k k'$ is isomorphic to the Koszul complex of $A'$ (which makes sense
since $k'$ is von Neumann regular). Thus $A'$ is Koszul.
\Edm

\setcounter{equation}{0}

\section{PBW theorem}
Our PBW theorem is the generalization to the $N$-case of a result (due to 
Braverman and Gaitsgory~\cite{bg:pbw}, see also Polishchuk and Positselski~\cite{pp:quad}) 
concerning non-homogeneous quadratic algebras over a field. It will be more convenient for us to
follow along the lines of~\cite{pp:quad}. 

Throughout this section, $k$ is a von Neumann regular ring, $V$ is a $k$-$k$-bimodule, and
$T(V)$ is the connected $\mathbb{N}$-graded $k$-$k$-algebra introduced in the previous section.
Actually, in this section, $T(V)$ will be seen as a connected $\mathbb{N}$-\emph{filtered}
$k$-$k$-algebra: the ring $T(V)$ is filtered by the sub-$k$-$k$-bimodules
$$F^n=\bigoplus_{0\leq i \leq n} V^{\otimes i},\ n \geq 0,$$
and the left or right $k$-actions coincide with the
products in the ring by elements of $F^0=k$. 

Our data are now an integer $N\geq 2$ and a sub-$k$-$k$-bimodule $P$ of $F^N$. Then the
two-sided ideal $I(P)$ is filtered by the $k$-$k$-bimodules 
$$I(P)^n=I(P)\cap F^n,\ n\geq 0.$$
We are interested in the
$\mathbb{N}$-filtered $k$-$k$-algebra $U=T(V)/I(P)$. The filtration of $U$ is formed by the
$k$-$k$-bimodules $U^n=F^n/I(P)^n$. Warning: $U^0$ can vanish (if $P=k$ for example), but the PBW
property will avoid this case.

Although $I(P)=\sum_{i,j\geq 0} V^{\otimes i}P\,V^{\otimes j}$, $I(P)^n$ may contain \emph{strictly}
the sum $\sum_{i+N+j\leq n} V^{\otimes i}P\,V^{\otimes j}$ in a non-trivial way. The following example is well-known.
\\ 

\addtocounter{Df}{1}
\Ea Take $n=N=2$, $k$ a field of characteristic 0, and $V$ finite-dimensional. Let
$f:V\times V
\rightarrow V$ be an alternate bilinear map, and let $P$ be the subspace of $F^2$ spanned by all
elements
$$r_{xy}= xy -yx- f(x,y), \ x,y \in V.$$
Set $c_{xyz}=[r_{xy},z]+[r_{yz},x]+[r_{zx},y]$, where $[\,,\,]$ is the usual commutator in $T(V)$.
Clearly, $c_{xyz}\in PV+VP$. Jacobi identity for $[\,,\,]$ shows that $c_{xyz}$ belongs to
$V^{\otimes 2}$, hence to $I(P)^2$. On the other hand,
$$c_{xyz}=-r_{f(x,y)z}-r_{f(y,z)x}-r_{f(z,x)y}-J(x,y,z),$$ 
where $J(x,y,z)=f(f(x,y),z)+f(f(y,z),x)+f(f(z,x),y) \in V$. Since $P\cap V=0$, $c_{xyz}$ belongs to
$P$ if and only if $J(x,y,z)=0$. Thus, if $(V,f)$ is not a Lie algebra (i.e., if $f$ does not
satisfy the Jacobi identity), $I(P)^2$ contains strictly $P$.
\\

The PBW theorem gives a ``simple'' (as the Jacobi identity) sufficient condition in order to have
\begin{equation} \label{but}
I(P)^n=\sum_{i+N+j\leq n} V^{\otimes i}P\,V^{\otimes j},\ \mathrm{for\ any}\ n\geq 0.
\end{equation}
By convention, the right-hand side of (\ref{but}) vanishes when $n<N$. So if (\ref{but})
holds, $I(P)^n=0$ (and $U^n= F^n$) when $n<N$. Our aim is now to express more conceptually
(\ref{but}). 

Introduce the $\mathbb{N}$-graded $k$-$k$-algebra $gr(U)$ associated to $U$, and denote by
$(gr(U)_n)$ its gradation. One has $gr(U)_n=U^n/U^{n-1}$. Using $I(P)^{n-1}=I(P)^n \cap
F^{n-1}$ and canonical isomorphisms, we identify $gr(U)_n$ to the $k$-$k$-bimodule $F^n/(I(P)^n
+F^{n-1})$. The product of the algebra $gr(U)$ is natural in these identifications. On the other
hand, let $\pi : F^N\rightarrow V^{\otimes N}$ be the projection associated to $F^N=
V^{\otimes N}\oplus F^{N-1}$. Then $R=\pi (P)$ is a sub-$k$-$k$-bimodule of $V^{\otimes N}$, and
we consider, as in the previous section, the connected $\mathbb{N}$-graded $k$-$k$-algebra
$A=T(V)/I(R)$. 

For each $n\geq 0$, let $\phi_n : V^{\otimes n} \rightarrow gr(U)_n$ be the composite of
canonical maps
$$V^{\otimes n} \rightarrow F^n \rightarrow F^n/(I(P)^n +F^{n-1}).$$
Since $P+F^{N-1}=R\oplus F^{N-1}$, one has $\phi_N(R)=0$. Thus the
surjective algebra morphism $\phi =
\oplus_{n\geq 0}\, \phi_n : T(V)\rightarrow gr(U)$ defines a surjective morphism of
$\mathbb{N}$-graded $k$-$k$-algebras
$$p: A \rightarrow gr(U).$$
Using the canonical
isomorphism
$$A_n=\frac{V^{\otimes n}}{I(R)_n} \cong \frac{V^{\otimes n}\oplus F^{n-1}}{I(R)_n\oplus
F^{n-1}}\ ,$$ 
one identifies $p_n :A_n \rightarrow gr(U)_n$ to the canonical map
$$\frac{F^{n}}{I(R)_n+ F^{n-1}} \rightarrow \frac{F^{n}}{I(P)^n+ F^{n-1}}\ .$$
Thus $p$ is an isomorphism if and only if $I(P)^n \subseteq I(R)_n + F^{n-1}$ for any $n\geq 0$.
The condition ``$p$ is an isomorphism'' is called the \emph{PBW property} (for $U$). 
\Bpo \label{car1}
The PBW property is equivalent to (\ref{but}) (actually, (\ref{but}) for $n\geq N-1$ suffices
since $I(P)^n\subseteq I(P)^{N-1}$ when $0\leq n \leq N-2$). 
\Epo
\Bdm
Introduce the notations
\begin{equation} \label{nota}
J^n =\sum_{i+N+j\leq n} V^{\otimes i}P\,V^{\otimes j},\ n\geq 0,
\end{equation}
with convention $J^n=0$ if $n<N$. One has $J^n \subseteq I(P)^n$. Assume firstly that $I(P)^n =J^n$
for any $n$. From that and $P\subseteq R+F^{N-1}$, one draws 
$$I(P)^n\subseteq \sum_{i+N+j=n} V^{\otimes i}R\,V^{\otimes j}+F^{n-1},$$
thus $p$ is an isomorphism. 

Conversely, assume that $p$ is an isomorphism. Prove equalities (\ref{but}) by induction on $n\geq
0$. Since $I(P)^0\subseteq I(R)_0+F^{-1}$, one has $I(P)^0=0$, hence (\ref{but}) for $n=0$. Let
$n$ be $\geq 1$, and assume (\ref{but}) for $n-1$. Since $I(P)^n \subseteq I(R)_n + F^{n-1}$ and
$R\subseteq P+F^{N-1}$, one has
$$I(P)^n \subseteq \sum_{i+N+j=n} V^{\otimes i}P\,V^{\otimes j} + F^{n-1},$$
where the sum is supposed to vanish if $n<N$ (because $I(R)_n=0$ in this case). Let $a$ be in
$I(P)^n$, and write down $a=b+c$, $b\in \sum_{i+N+j=n} V^{\otimes i}P\,V^{\otimes j}$, $c\in
F^{n-1}$. Clearly,
$c=a-b$ belongs to $I(P)\cap F^{n-1}=I(P)^{n-1}$, hence to $J^{n-1}\subseteq J^n$ by induction
hypothesis. But $b\in J^n$, thus $a\in J^n$. 
\Edm
\Bpo \label{car2}
Keeping notations (\ref{nota}), the PBW property is equivalent to 
\begin{equation} \label{vbut}
J^n \cap F^{n-1} =J^{n-1},\ \mathrm{for\ any}\ n\geq N.
\end{equation}
\Epo
\Bdm
If $p$ is an isomorphism, (\ref{but}) shows that $J^n\cap F^{n-1}= I(P)^n \cap F^{n-1}=
I(P)^{n-1}= J^{n-1}$. Conversely, assume (\ref{vbut}). Fix $n\geq N-1$. Since $I(P)= \sum_{i,j\geq
0} V^{\otimes i}P\,V^{\otimes j}$, one has
$$I(P)^n= ( \bigcup_{i\geq N}J^i) \cap F^n=\bigcup_{i\geq N}(J^i\cap F^n). $$ 
For $N\leq i \leq n$, $J^i \subseteq J^n$ hence $J^i \cap F^n \subseteq J^n$. Next
$J^{n+1}\cap F^n=J^n$ because of (\ref{vbut}) and $n+1\geq N$. Moreover $J^{n+2}\cap F^n=
J^{n+2}\cap F^{n+1}\cap F^n = J^{n+1}\cap F^n= J^n$. A straightforward induction shows that $J^i\cap
F^n =J^n$ for any $i\geq n+1$. Thus $I(P)^n\subseteq J^n$, and we conclude by the previous
proposition. 
\Edm
\\

For $n=N$ and $n=N+1$, (\ref{vbut}) is respectively equivalent to
\begin{equation} \label{I}
P\cap F^{N-1} =0,
\end{equation}
\begin{equation} \label{J}
(PV+VP) \cap F^{N} \subseteq P.
\end{equation} 
(For $n=N+1$, one can replace $P+PV+VP$ by $PV+VP$ 
in the non-trivial inclusion of (\ref{vbut}).) 

When $N=2$ and $k$ is a field,
(\ref{I}) and (\ref{J}) are the conditions (I) and (J) of Braverman and
Gaitsgory. Their 
PBW theorem
asserts that if $A$ is Koszul, (I) and (J) are sufficient in order to have the PBW property. This
theorem extends to any $N\geq 2$ and any von Neumann regular ring $k$, as follows.
\Bte \label{PBW}
\emph{(PBW theorem in the $N$-case)} Assume that $k$ is a von Neumann regular ring, $V$ is a
$k$-$k$-bimodule, $N$ is an integer $\geq 2$, and $P$ is a sub-$k$-$k$-bimodule of $F^N$, where
$F^n=\oplus_{0\leq i \leq n} V^{\otimes i}$ for any $n \geq 0$. Set $U=T(V)/I(P)$ and
$A=T(V)/I(R)$, where $R=\pi(P)$ and $\pi$ is the projection of $F^N$ onto $V^{\otimes N}$ modulo
$F^{N-1}$. 

Assume that the graded left $k$-module Tor$_3^A(k_A,_Ak)$ is concentrated in degree $N+1$.
Then (\ref{I}) and (\ref{J}) imply that the PBW property holds, i.e., the natural algebra morphism
$p:A\rightarrow gr(U)$ is an isomorphism.
\Ete

Before proving the theorem, let us introduce the arrows $\varphi^{i,i+N-1}$ which are a generalization
of the arrows $\psi^{i,i+1}$ of~\cite{pp:quad}. From now on, we assume (\ref{I}). As
$\ker(\pi)=F^{N-1}$, $\pi$ realizes an isomorphism from $P$ onto $R$. For any $x\in R$, let $y$ be
the element of $P$ such that $\pi(y)=x$, and set
$y=x-\varphi(x)$. That defines a $k$-$k$-linear map
$\varphi : R\rightarrow F^{N-1}$. Then $P$ can be described only by $R$ and $\varphi$:
\begin{equation} \label{psi}
P=\{ x-\varphi(x);\ x\in R\}.
\end{equation}
Fix $i\geq 1$, $j\geq 0$. Using the identifications introduced in Section 2 just before
Lemma \ref{formulas}, the $k$-$k$-linear map 
$$1_{V^{\otimes(i-1)}}\otimes_k \varphi \otimes_k 1_{V^{\otimes j}} : V^{\otimes(i-1)} \otimes_k
R\otimes_k V^{\otimes j} \rightarrow V^{\otimes(i-1)} \otimes_k F^{N-1} \otimes_k V^{\otimes j}$$
is identified to a $k$-$k$-linear map $V^{\otimes(i-1)}  R V^{\otimes j} \rightarrow T(V)$
which is denoted by $\varphi^{i,i+N-1}$ (actually it arrives in $V^{\otimes(i-1)}F^{N-1}V^{\otimes j}$
which is naturally embedded in $T(V)$). 
\Bpo \label{car3}
Assume that (\ref{I}) holds. Then (\ref{J}) is equivalent to 
\begin{equation} \label{J'}
(\varphi^{1,N}-\varphi^{2,N+1})(W_{N+1}) \subseteq P.
\end{equation}
\Epo
\Bdm
Let $a$ in $W_{N+1}$. Then $a-\varphi^{1,N}(a) \in PV$ and $a-\varphi^{2,N+1}(a) \in VP$, so
$(\varphi^{1,N}-\varphi^{2,N+1})(a)$ belongs to $(PV+VP) \cap F^{N}$. Conversely, if $x$ belongs to
$(PV+VP) \cap F^{N}$, decompose
$$x=\sum_i (x_i-\varphi(x_i))v_i + \sum_j v'_j (x'_j-\varphi(x'_j)),$$
with $x_i$, $x'_j$ in $R$, $v_i$, $v'_j$ in $V$. Since $x\in F^N$, $\sum_i
x_i v_i + \sum_j v'_j x'_j =0$, hence
$$x=-(\varphi^{1,N}-\varphi^{2,N+1})(\sum_i x_i v_i),$$
where $\sum_i x_i v_i \in W_{N+1}$. Thus $(\varphi^{1,N}-\varphi^{2,N+1})(W_{N+1})= (PV+VP) \cap F^{N}$.
\Edm
\\ \\
\emph{Proof of Theorem \ref{PBW}.} According to Proposition \ref{car2}, it suffices to prove that
$J^n \cap F^{n-1} \subseteq J^{n-1}$ by induction on $n\geq N$. Assume $n\geq N+2$ and the inclusion
true for $n-1$. Let $x$ be in $J^n \cap F^{n-1}$. In the sequel, the symbol $\equiv$ means equality
modulo $J^{n-1}$. One can find $x_i$ in $V^{\otimes (i-1)}RV^{\otimes (n-N-i+1)}$ for $1\leq i \leq
n-N+1$ such that
$$x \equiv \sum_{i=1}^{n-N+1} (x_i-\varphi^{i,i+N-1}(x_i)).$$
Since $x\in F^{n-1}$, the sum of the $x_i$'s vanishes, so 
$$x_{n-N+1}\in (V^{\otimes (n-N)}R)\cap (I(R)_{n-1}V).$$
But Tor$_3^A(k_A,_Ak)$ is concentrated in degree $N+1$. Theorem \ref{tor3} shows that 
there exist $y_i$ in $V^{\otimes (i-1)}RV^{\otimes (n-2N-i+1)}R$ for $1\leq i \leq n-2N+1$ and
$y_{n-N}$ in $V^{\otimes (n-N-1)}W_{N+1}$ such that 
$$x_{n-N+1}= \sum_{i=1}^{n-2N+1} y_i\ +y_{n-N}.$$
Note that if $n<2N$, the latter equality reduces to $x_{n-N+1}=y_{n-N}$, so that the other $y_i$'s
are considered as vanishing in this case (actually, (ec) and (\ref{other}) can be put together
in the statement of Theorem \ref{tor3}). Using (\ref{J'}), we get
$$\varphi^{n-N+1,n}(x_{n-N+1})\equiv \sum_{i=1}^{n-2N+1} \varphi^{n-N+1,n}(y_i)\ +
\varphi^{n-N,n-1}(y_{n-N}).$$

Fix the index $i$, $1\leq i \leq n-2N+1$. Decompose
$$\varphi^{i,i+N-1} : V^{\otimes (i-1)}RV^{\otimes (n-2N-i+1)}R \rightarrow V^{\otimes
(i-1)}F^{N-1}V^{\otimes (n-2N-i+1)}R$$
as $\varphi^{i,i+N-1} = \sum_{j=0}^{N-1} \varphi_j^{i,i+N-1}$, where
$$\varphi_j^{i,i+N-1} : V^{\otimes (i-1)}RV^{\otimes (n-2N-i+1)}R \rightarrow V^{\otimes (n-2N+j)}R.$$
But the two $R$'s in the tensor product $E=V^{\otimes (i-1)}RV^{\otimes (n-2N-i+1)}R$ do not
overlap. Thus the equality
$$\varphi^{i,i+N-1}\varphi^{n-N+1,n} =  \sum_{j=0}^{N-1}\varphi^{n-2N+j+1,n-N+j} \varphi_j^{i,i+N-1}$$ 
holds \emph{on E}, so that \emph{on E}, one has 
$$\varphi^{i,i+N-1}-\varphi^{n-N+1,n} =  \sum_{j=0}^{N-1}(\mathrm{id}-\varphi^{n-2N+j+1,n-N+j})
\varphi_j^{i,i+N-1}-(\mathrm{id}-\varphi^{i,i+N-1})\varphi^{n-N+1,n}$$
and we get $\varphi^{n-N+1,n}(y_i) \equiv \varphi^{i,i+N-1}(y_i)$. Finally
$$x\equiv -\sum_{i=1}^{n-2N+1} \varphi^{i,i+N-1}(x_i+y_i) -
\sum_{i=n-2N+2}^{n-N-1}\varphi^{i,i+N-1}(x_i)  - \varphi^{n-N,n-1}(x_{n-N}+ y_{n-N})$$
in which the second sum does not occur if $N=2$. 

Set $x'_i=x_i+y_i$ for $1\leq i \leq n-2N+1$ and $i=n-N$. Set $x'_i=x_i$ for $n-2N+2 \leq i \leq
n-N-1$. Then for $1\leq i \leq n-N$, $x'_i$ belongs to $V^{\otimes (i-1)}RV^{\otimes (n-N-i)}V$,
and $\sum_{i=1}^{n-N} x'_i =0$. Introducing
$$x'= \sum_{i=1}^{n-N} (x'_i-\varphi^{i,i+N-1}(x'_i)),$$
we get $x\equiv x'$ and $x' \in (J^{n-1}V)\cap (F^{n-2}V)$. Lemma \ref{formulas} gives
$(J^{n-1}V)\cap (F^{n-2}V)=(J^{n-1}\cap F^{n-2})V$, therefore $x'$ belongs to $J^{n-2}V$ by
induction hypothesis. Thus $x'$ and $x$ are in $J^{n-1}$. 
\Edm
\Bpo \label{car4}
Assume that (\ref{I}) holds. Let $\varphi :R \rightarrow F^{N-1}$ be the $k$-$k$-linear map such
that $P=\{x-\varphi(x);\ x \in R \}$. Decompose $\varphi = \sum_{j=0}^{N-1}\varphi_j$, $\varphi_j : R\rightarrow
V^{\otimes j}$. Then (\ref{J}) is equivalent to all following relations
\begin{equation} \label{J'1}
(\varphi^{1,N}_{N-1}-\varphi^{2,N+1}_{N-1})(W_{N+1})\subseteq R,
\end{equation}
\begin{equation} \label{J'2}
\left(\varphi_j(\varphi^{1,N}_{N-1}-\varphi^{2,N+1}_{N-1})+\varphi^{1,N}_{j-1}-\varphi^{2,N+1}_{j-1}
\right)(W_{N+1})=0,\ 1\leq j \leq N-1,
\end{equation}
\begin{equation} \label{J'3}
\varphi_0(\varphi^{1,N}_{N-1}-\varphi^{2,N+1}_{N-1})(W_{N+1})=0.
\end{equation}
\Epo
\Bdm
Apply Proposition \ref{car3}. Let $x$ be in $W_{N+1}$ and $X=(\varphi^{1,N}-\varphi^{2,N+1})(x)$. Since the
two projections in $P\subseteq R\oplus F^{N-1}$ are respectively $\pi$ and $-\varphi \pi$, then $X$
belongs to $P$ if and only if $\pi(X)\in R$ and $X=\pi (X)-\varphi \pi(X)$. The component of degree 
$j$, $0\leq j \leq N$, of $X$ is $(\varphi^{1,N}_{j-1}-\varphi^{2,N+1}_{j-1})(x)$ (vanishing if
$j=0$). In particular  
$$\pi (X)= (\varphi^{1,N}_{N-1}-\varphi^{2,N+1}_{N-1})(x).$$
And the component of degree $j\leq N-1$ of $\varphi \pi(X)$ is
$\varphi_j(\varphi^{1,N}_{N-1}-\varphi^{2,N+1}_{N-1})(x)$.
\Edm
\\ 

\addtocounter{Df}{1}
\Ea In Example 3.1, (\ref{I}) holds and $\varphi_0 =0$, so (\ref{J'3}) holds. As elements of
$W_3$ are totally skew-symmetric, (\ref{J'1}) is easily checked, and (\ref{J'2}) is equivalent to
the Jacobi identity for $f$.
\\ 

\addtocounter{Df}{1}
\Ea (Down-up algebras~\cite{br:dup}) Here $k$ is a field, $\alpha$, $\beta$, $\gamma$ are in $k$
with $\beta\neq 0$, $U=U(\alpha, \beta, \gamma)$ is the associative $k$-algebra with two generators
$d$, $u$ and following relations
$$d^2u = \alpha dud + \beta ud^2 + \gamma d,$$
$$du^2 = \alpha udu + \beta u^2d + \gamma u.$$
Taking $d$ and $u$ of degree 1, $U$ is a non-homogeneous cubic algebra, and the homogeneous
cubic algebra $A$ is defined by relations $r_1=r_2=0$, where $r_1=d^2u- \alpha dud - \beta ud^2$
and $r_2=du^2 - \alpha udu - \beta u^2d$. Then $A$ is the AS-regular algebra of global dimension 3
which is cubic of type S$_1$~\cite{as:regular}. In particular, $A$ is Koszul~\cite{bm:dual}.
Moreover~\cite{bm:dual}, $W_4= RV\cap VR$ is one-dimensional and a generator is
$$w=r_1 u-\beta r_2d = -\beta ur_1+d r_2.$$ 
Let us show the PBW property. Firstly, $P\cap F^2=0$ is clear. Use Proposition
\ref{car4}. One has $\varphi_0=\varphi_2=0$ and $\varphi_1$ is defined by $\varphi_1(r_1) =\gamma d$,
$\varphi_1(r_2) =\gamma u$. Thus (\ref{J'1}) and (\ref{J'3}) hold, whereas (\ref{J'2}) comes from the
calculation
$$(\varphi^{1,3}_{1}-\varphi^{2,4}_{1})(w)=\gamma d u -\beta \gamma ud-( -\beta u \gamma d + d\gamma
u)=0.$$
The fact that $p: A \rightarrow gr(U)$ is an isomorphism allows us to deduce properties of $U$ from
those of $A$. For example, using~\cite{atv:grot, atv:mod}, $U(\alpha,
\beta, \gamma)$ is a noetherian domain of Gelfand-Kirillov dimension 3
(see also~\cite{br:dup}). 
\Bdf \label{koszul}
Assume that $k$ is a von Neumann regular ring, $V$ is a
$k$-$k$-bimodule, $N$ is an integer $\geq 2$, and $P$ is a sub-$k$-$k$-bimodule of $F^N$, where
$F^n=\oplus_{0\leq i \leq n} V^{\otimes i}$ for any $n \geq 0$. Set $U=T(V)/I(P)$ and
$A=T(V)/I(R)$, where $R=\pi(P)$ and $\pi$ is the projection of $F^N$ onto $V^{\otimes N}$ modulo
$F^{N-1}$. Then $U$ is said to be Koszul if the graded algebra (with
$N$-homogeneous relations) $A$ is Koszul and if the PBW property
holds.
\Edf
Note that if $U$ is Koszul, the integer $N$ is uniquely determined,
and we shall say that $U$ is $N$-Koszul. Definition \ref{koszul} generalizes the
definition of the graded situation: if the relations of $U$ are all
$N$-homogeneous, then $U=A$ is graded and the PBW property holds
trivially. 
\\ 

\addtocounter{Df}{1}
\Rm \label{trivial} Let us show how $\varphi_0$ measures the obstruction for $k$ to be an $U$-module. Assume firstly that $\varphi_0=0$. It is easy to check that the natural projection $q :T(V)\rightarrow k$ vanishes on the ideal $I(P)$, so that we get a natural $k$-$k$-linear map $\eps_U : U\rightarrow k$. Then the left $U$-module $k$ deduced from $\eps_U$ satisfies the following properties :

(1) $u.1=0$ for any $u\in V^{\o n}$ and $1\leq n \leq N-1$,

(2) $\bar{x}.1=0$ for any $x\in R$, where $\bar{x}$ denotes the class of $x$ in $U^N$.
\\
Conversely, assume that $k$ is a $U$-module satisfying (1) and (2). Then it is immediate that $\varphi_0=0$.

\setcounter{equation}{0}

\section{$N$-inhomogeneous
  algebras associated to finite groups}\label{gp}
Throughout this section, $k$ is an algebraically closed field of characteristic 0,
$V$ a $k$-vector space of dimension $n$, and
$\Gamma$ is a finite subgroup of $GL(V)$. The group algebra $K=k[\Gamma]$
is a semi-simple ring, hence is von Neumann regular. Fix $p$, $1< p
\leq n$. 

Given a $k$-linear map $\psi : \Lambda_k ^{p}(V) \rightarrow K$,
put 
$$H_{\psi} = (T_k(V)\# \Gamma) / I(\mathrm{Alt}(v_1,\ldots,
v_{p})-\psi (v_1,\ldots, v_{p})\, ; \ v_1,\ldots, v_{p} \in V),$$
where Alt stands for the anti-symmetrization in $T_k^{p} (V)$. Consider
the $K$-$K$-bimodule $E=V\otimes_k K$, with left $\Gamma$-action given
by $g : v\otimes a \mapsto g(v)\otimes (ga)$, and right
$\Gamma$-action given by $v\otimes a \mapsto v\otimes
(ag)$,
where $ga$ and $ag$ stand for the product in the group algebra. Then
$T_K(E) \cong T_k(V)\# \Gamma$, so that $H_{\psi}$ is identified to
the $\mathbb{N}$-filtered $K$-$K$-algebra $T_K(E)/I(P)$
where $P$ is the sub-$K$-$K$-bimodule of $T_K(E)$ generated by
the following elements of $T_K^{p} (E)\oplus K$:
$$\mathrm{Alt}(v_1,\ldots, v_{p})-\psi (v_1,\ldots, v_{p}) \
\mathrm{with}  \ v_1,\ldots, v_{p} \in V.$$
So $H_{\psi}$ has relations of degree $\leq N=p$. Using notations
of the previous section, $U=H_{\psi}$ and $A=H_{\psi=0}$. Clearly, $A$ is obtained from
the $N$-homogeneous $k$-$k$-algebra
$$\mathcal{A}=T_k(V) / I(\mathrm{Alt}(v_1,\ldots, v_{p})\, ; \ v_1,\ldots, v_{p}
\in V)$$ by the natural change of rings $k \rightarrow K$. It is known
that $\mathcal{A}$ is Koszul
\cite{rb:nonquad}.
Thus Proposition \ref{change} shows that $A$ is Koszul.

\Blm In the notation of \eqref{J'1}-\eqref{J'3}, in $T_K^{p+1}E$,
we have
$W_{p+1}=(\wedge^{p+1}_kV)\otimes_k K$ as a $K$-subbimodule in
$T_k^{p+1}V\otimes_k K$.
\Elm

\Bdm According to  \cite{rb:nonquad}, the algebra $\mathcal{A}$
is Koszul. Furthermore, if $\mathcal{R}$ denotes the
space of relations of $\mathcal{A}$, then the space of relations of
$A$ equals $R=\mathcal{R}\otimes_kK$. Hence we get
$W_q=\cap_{i+p+j=q}E^{\o i}\o R\o E^{\o j}= (\wedge^q_kV)\o_kK,$
as a $K$-subbimodule in $T^q_kV\o_kK.$\qed
\medskip

Now, we may write the map $\psi : \Lambda_k ^{p}(V) \rightarrow K$ 
in the form $\psi=\sum_{g\in\G}\psi_g\cdot g,$ where
$\psi_g: \Lambda_k ^{p}(V) \rightarrow k$ are certain
linear maps.
Further, let  $\Gamma$ act on $K$ by conjugation.

\Blm\label{lem1} The algebra $H_{\psi}$ is Koszul
if and only if $\psi : \Lambda_k ^{p}(V) 
\rightarrow K$ is $\Gamma$-equivariant
and, for any $g\in \G$ and  $v_1,\ldots, v_{p+1}\in V,$
in $V$ one has
\beq{identity}
\sum_{i=1}^{p+1}
(-1)^i\cdot  \psi_g(v_1,\ldots,v_{i-1},v_{i+1},\ldots, v_{p+1})\cdot
(\Id-
(-1)^pg)(v_i)=0.
\eeq
\Elm
\Bdm
First of all, we
claim that if the PBW-property holds then
the map $\psi $ must be $\G$-equivariant.
To see this, let  $e_\rho\in k[\G]$ denote the
central idempotent corresponding to an
 irreducible representation
$\rho$ of $\G$. Let $\G$ act on $K$ via the adjoint action and
on the vector space $\Hom_k(T^p_kV,T^p_kV\o_k K)$ via the corresponding induced action.
The map $\mathrm{Alt}: v_1,\ldots,v_p\mapsto \mathrm{Alt}(v_1,\ldots,v_p)\o1$
is clearly $\G$-equivariant, hence, in  $\Hom_k(T^p_kV,T^p_kV\o_k K)$,
we have $e_\rho(\mathrm{Alt})=0$, for every {\em nontrivial}
 irreducible representation
$\rho$ of $\G$. Now, if the map $\psi$ is not
 $\G$-equivariant, then there exists a  nontrivial irreducible representation
$\rho$ such that $e_\rho(\psi)\neq 0$.
We conclude that there exist $v_1,\ldots,v_p\in V$
 such that $\bigl(e_\rho(\mathrm{Alt}-\psi)\bigr)(v_1,\ldots,v_p)=
\bigl(e_\rho(\psi)\bigr)(v_1,\ldots,v_p)\neq 0$.
This means that in the algebra $H_\psi$
we have a relation $a=0$, where
$a:=\bigl(e_\rho(\psi)\bigr)(v_1,\ldots,v_p)\in K$
is a nonzero element. Thus, the canonical map
$K\to H_\psi$ is not injective and PBW-property fails. Actually, denoting by $(F^n)_{n\geq 0}$ the natural filtration of $T_K(E)$, it is elementary that $P\cap F^{p-1}=0$ (i.e., the condition (3.4)) holds if and only if $\psi : \Lambda_k ^{p}(V) 
\rightarrow K$ is $\Gamma$-equivariant.

The rest of the argument is very similar to the proof of
\cite[formula (2.3)]{eg:reflec}. Assume that $\psi$ is $\Gamma$-equivariant.
Since $\varphi=\varphi_0=\psi$,
the criterion of  Proposition \ref{car4} (i.e., the condition (3.5)) reduces to the
following equation in $V$:
$$(\psi^{1,p}-\psi^{2,p+1})(\mathrm{Alt}(v_1,\ldots,v_p))=0,
\quad\forall 
v_1,\ldots,v_{p+1}\in V.$$
Explicitly, writing $\eps(\si)$ for the sign of permutation
$\si$, the equation reads
$$
\sum_{\sigma\in \Si_{p+1}}
\eps(\si)\left[\psi(v_{\si(1)},\ldots,v_{\si(p)})\cdot
v_{\si(p+1)}-
v_{\si(1)}\cdot\psi(v_{\si(2)},\ldots,v_{\si(p+1)})\right]=0.
$$
Rewriting this expression one obtains the following condition:
\begin{align*}
0&=\sum_{i=1}^{p+1}(-1)^{p+1-i}
\sum_{\tau\in \Si_p}
\eps(\tau)\psi(v_{\tau(1)},\ldots,v_{\tau(i-1)},
v_{\tau(i+1)},\ldots,v_{\tau(p+1)})\cdot
v_i\\
&-\sum_{i=1}^{p+1}(-1)^{i-1}\sum_{\tau\in \Si_p}
\eps(\tau)v_i\cdot\psi(v_{\tau(1)},\ldots,v_{\tau(i-1)},
v_{\tau(i+1)},\ldots,v_{\tau(p+1)})\\
&=
\sum_{i=1}^{p+1}(-1)^i
[\psi(v_{\tau(1)},\ldots,v_{\tau(i-1)},
v_{\tau(i+1)},\ldots,v_{\tau(p+1)}),\,v_i]_\pm,
\end{align*}
where we use the notation 
$[a,v]_\pm=a\cdot v-(-1)^p v\cdot a,$
for any $a\in K,\, v\in V$.

Further, for any $v\in V$ and $g\in \G$,
in $T_k(V)\# \Gamma$, we have 
$$[v,g]_{\pm}= (v\o 1)\cdot g -(-1)^p g\cdot (v\o 1)= v\o g - (-1)^p g(v)\o g= (v- (-1)^p g(v))\o g.$$
Therefore, writing
$\psi=\sum_{g\in\G}\psi_g\cdot g,$
the last displayed formula reads
\begin{align*}
0&=\sum_{i=1}^{p+1}(-1)^i\cdot \Bigl(
\sum\nolimits_{g\in\G}
\psi_g(v_1,\ldots,v_{i-1},v_{i+1},\ldots, v_{p+1})\cdot (v_i-(-1)^p g(v_i))\cdot g\Bigr)\\
&=\sum_{g\in\G}\Bigl(\sum_{i=1}^{p+1}(-1)^i\cdot 
\psi_g(v_1,\ldots,v_{i-1},v_{i+1},\ldots, v_{p+1})\cdot (\Id-(-1)^pg)(v_i)\Bigr)\cdot g.
\end{align*}
Thus, the coefficient in front of each element $g$ in the second line of
the formula must vanish, and the Lemma follows.
\qed

To simplify notation, we write
$\otimes=\otimes_k$ and $\Lambda^i(-)=\Lambda_k ^i(-),$
etc.

Next, fix $g\in\G$.
Write  $M=M_g:=\op{Image}(\Id-(-1)^p g),$ and $L=L_g:=\op{Ker}(\Id-(-1)^p g),$
and put $m:=\dim M,\, l:=\dim L$.
Since $g$ is an element of $GL(V)$ of finite order,
hence, semisimple, we have a direct sum decomposition
$V=M\oplus L$. 
Thus, we have $m+l=n$ and 
we have a natural isomorphism
\beq{tensor}
\Lambda ^{p}(V)
=\Lambda ^{p}(M\oplus L)\cong\oplus_{i+j=p}
\Lambda ^i(M)\otimes\Lambda ^j(L).
\eeq
Note also that the space $\Lambda ^m(M)$
is 1-dimensional.

\Blm\label{lem2} A linear map $\phi:\Lambda ^{p}(V)=\Lambda ^{p}(M\oplus L) \rightarrow k$ 
satisfies the equation

$$\sum_{i=1}^{p+1}
(-1)^i\cdot  \phi(v_1,\ldots,v_{i-1},v_{i+1},\ldots, v_{p+1})\cdot (\Id-(-1)^p g)(v_i)=0,
\quad\forall 
v_1,\ldots, v_{p+1}\in V.
$$
if and only if $\phi=\phi^{m,p-m}$, where $\phi^{m,p-m}$ is the restriction of $\phi$ to $\Lambda ^m(M)\otimes\Lambda ^{p-m}(L)$. In particular, if $\phi$ does not vanish, then $p\geq m$.
\Elm
\Bdm
Recall that, for any $d$-dimensional vector
space $E$ and an integer $0\leq p\leq d$, one has a  {\em Koszul 
exact sequence} (a version of Koszul complex):
\begin{align}\label{kos}
K^p(E):\;0\map\Lambda^p(E^*)
\stackrel{\partial^p_E}\map&
\Lambda^{p+1}(E^*)\otimes E
\map
\Lambda^{p+2}(E^*)\otimes S^2(E)\map\\
&\ldots\map
\Lambda^d(E^*)\otimes S^{d-p}(E)\map k\map 0.\nonumber
\end{align}
The first differential
$\partial^p_E: \Lambda^p(E^*)
\to\Lambda^{p+1}(E^*)\otimes E=\Hom(\Lambda^{p+1}E^*, E)$
is defined
by the formula
$$(\partial^p_E\alpha)(v_1,\ldots, v_{p+1}):=
\sum_{i=1}^{p+1}
(-1)^i\cdot  \alpha(v_1,\ldots,v_{i-1},v_{i+1},\ldots, v_{p+1})\cdot
v_i,
$$
for any
$v_1,\ldots, v_{p+1}\in E.$
The other differentials are given by a similar
formula. Thus, from \eqref{kos} we see that
\beq{injective}
\text{The map $\partial^p_E$ is injective unless
$p=\dim E$.}
\eeq

The Koszul exact sequence \eqref{kos} is compatible with
tensor products. In particular,
for $V=M\oplus L$
we have $\partial^p_V=\partial_M\otimes\Id_L\pm \Id_M\otimes\partial_L$.
In more detail,  using the
direct sum decomposition
\eqref{tensor}, the differential
$$\partial^p_V:\ 
\Lambda^p(M^*\oplus L^*)\to
\Lambda^{p+1}(M^*\oplus L^*)\otimes (M\oplus L)$$
may be written as a sum
$\partial^p_V=\sum_{r+s=p}\partial_V^{r,s}$,
of the following components
$$\partial_V^{r,s}:\ 
\Lambda^r(M^*)\otimes\Lambda^s(L^*)
\map \oplus_{i+j=p+1}
\Lambda ^i(M)\otimes\Lambda ^j(L)\otimes (M\oplus L).
$$
With this notation, for any $r,s$, one has a Leibniz type formula:
\beq{rs}
\partial_V^{r,s}=\partial_M^r\otimes\Id^s_L +
(-1)^r\cdot\Id^r_M\otimes\partial_L^s,
\eeq
where the two summunds are interpreted as maps
$$\partial_M^r\otimes\Id^s_L:\
\Lambda^r(M^*)\otimes\Lambda^s(L^*)\map
\Lambda^{r+1}(M^*)\otimes\Lambda^s(L^*)\otimes M,
$$
resp.,
$$
\Id^r_M\otimes\partial_L^s
:\
\Lambda^r(M^*)\otimes\Lambda^s(L^*)\map
\Lambda^r(M^*)\otimes\Lambda^{s+1}(L^*)\otimes L.
$$
The target spaces in these formulas are viewed
as being subspaces of the vector space
$\oplus_{i+j=p+1}
\Lambda ^i(M)\otimes\Lambda ^j(L)\otimes (M\oplus L).$

We can now complete the proof of the Lemma.
Observe that, in terms of differential
$\partial_V^p$, equation
in the statement of the Lemma reads
$(\Id-(-1)^p g)(\partial_V^p(\phi))=0$.
Now, the map $(\Id-(-1)^p g)$ takes the subspace $L\sset L\oplus M=V$ to zero,
and restricts to an invertible operator
$M\to M$.
Therefore, writing
$$\pi_M:
\oplus_{i+j=p+1}
\Lambda ^i(M)\otimes\Lambda ^j(L)\otimes (M\oplus L)\map
\oplus_{i+j=p+1}
\Lambda ^i(M)\otimes\Lambda ^j(L)\otimes M
$$
for the natural projection,
we obtain
$(\Id-(-1)^p g)(\partial_V^p(\phi))=0
\enspace\Longleftrightarrow\enspace
\pi_M(\partial_V^p(\phi))=0$.
Further, it is clear that we have
$$\pi_M\circ (\partial_M^r\otimes\Id^s_L)=
(\partial_M^r\otimes\Id^s_L),\quad
\text{and}\quad
\pi_M\circ (\Id^r_M\otimes\partial_L^s)=0.
$$
Thus, using formula
\eqref{rs}, we conclude that
\beq{van}
(\Id-(-1)^p g)(\partial_V^p(\phi))=0
\enspace\Longleftrightarrow\enspace
(\partial^r_M\otimes\Id_L^{p-r})(\phi^{r,p-r})=0,
\enspace\forall 0\leq r\leq p.
\eeq

Next, observe that the map $\partial^r_M$,
hence the map $(\partial^r_M\otimes\Id_L^{p-r})$,
is injective
for all $r<\dim M$, by \eqref{injective}.
Therefore, the vanishing condition on the
right of \eqref{van} holds if and only
if $r=\dim M=m$. Furthermore, in the latter case,
the space $\Lambda^rM$ is 1-dimensional.
It follows that the map $\phi^{m,p-m}$ may be written
as a tensor product $\phi_M\otimes \phi_L$ of linear maps
$\phi_M: \Lambda^mM\to k$ and $\phi_L: \Lambda^{p-m}L\to k$.
The Lemma is proved.
\qed

Now, for any $g\in \G$ we put $a(g):=\dim M_g$. Recall that $M_g= \op{Im}(\Id-(-1)^p g)$ and $L_g=\op{Ker}(\Id-(-1)^p g)$.
Lemmas \ref{lem1} and \ref{lem2} yield the following result.

\Bte \label{maint}
The $\mathbb{N}$-filtered $K$-$K$-algebra $H_{\psi}$ is Koszul
if and only if $\psi =\sum_{g\in\G}\psi_g\cdot g,$ : $\Lambda_k ^{p}(V) \rightarrow K$ is $\Gamma$-equivariant and for any $g\in \Gamma$, $\psi_g$ vanishes in all components $\Lambda ^i(M_g)\otimes\Lambda ^{p-i}(L_g)$ of $\Lambda ^{p}(V)$ unless if $i=a(g)$.

In particular, $H_{\psi}$ is Koszul if $\psi$
has the form $\psi=\sum_{g\in \G}
(\psi'_g\otimes\psi''_g)\cdot g,$ 
where
\begin{align*}
&\psi':\ g\mapsto \psi'_g\in \Lambda^{a(g)}\bigl(M_g\bigr)^{\ast},
\quad\text{resp.,}\\
&
\psi'':\ g\mapsto \psi''_g\in
\Lambda^{p-a(g)}\bigl(L_g\bigr)^*,
\end{align*}
are some $\G$-equivariant maps
(where $\Gamma$ acts on itself by conjugation). \qed
\Ete

The following special case of the Theorem arises in
many interesting examples.
Let $\phi : \Lambda^pV\to k$ be a $\G$-invariant
linear function.
For each $g\in \G$ let
\begin{align*}
&\psi_g = \phi  \quad \mathrm{in} \quad \Lambda ^{a(g)}(M_g)\otimes\Lambda ^{p-a(g)}(L_g)\\
&\psi_g = 0     \quad \mathrm{in\ any\ other\ component} \quad \Lambda ^i(M_g)\otimes\Lambda ^{p-i}(L_g).
\end{align*}
Note that, in general, $\psi_g \neq \phi$ 
unless $M_g= V$ or 0.

Now put $\psi=\sum_{g\in \G}
\psi_g \cdot g\in \Hom_{\G}(\Lambda^pV,K)$.

\Bcr \label{mainc} With the above notations,
the algebra $H_\psi$ is Koszul.
\Ecr

Symplectic reflection algebras appear as a special case
of this construction where $p=2$ and $\phi$ is a symplectic
2-form on $V$. Notice that the previous corollary still holds if $\psi$ is multiplied by any map $m : \G \rightarrow \mathbb{C}$ which is constant on any conjugation class.

\setcounter{equation}{0}

\section{Koszul complex}\label{kc}
Notations and assumptions are those of Definition \ref{koszul}. So $U$ is a filtered $N$-Koszul algebra over the von Neumann regular ring $k$. We are interested in defining the Koszul complex of $U$ for bimodules. Braverman and Gaitsgory defined it in the case $N=2$ and $\varphi=\varphi_1$ as a subcomplex of the bar resolution (see 5.4 in \cite{bg:pbw}). But it seems hard to proceed in the same manner when $N>2$, because $U\o R \o U$ is not naturally included in $U\o U^{\o2} \o U$. 

Having in mind the situation of Section \ref{gp}, we limit ourselves to the case $N\geq 2$ and $\varphi=\varphi_0$. Our method is to construct an $N$-differential on $(U\o W_n\o U)_{n\geq 0}$, and then to make an adequate contraction in order to get the Koszul complex as in \cite{bdvw:homog,bm:dual}. Recall the meaning of the notation $W_n$ for any $n\geq 0$:
$$W_n=\bigcap_{i+N+j=n}V^{\otimes i}RV^{\otimes j}.$$

Denote by $\mu : U\o U \rightarrow U$ the multiplication of $U$, and denote by $\mu|_{U\o V}$ its restriction to $U\o V$. Define the $U$-$k$ linear map
$$d_l : U\o W_n \rightarrow U\o V^{\o (n-1)}$$
as the restriction of $\mu|_{U\o V} \o 1_{V^{\o(n-1)}}$ to $U\o W_n$. Actually, $d_l$ maps into $U\o W_{n-1}$. In fact, for any $w\in W_n$, one can write down $w=\sum_j v_jw_j$ where $v_j \in V$ and $w_j\in W_{n-1}$, so that $d_l(a\o w)=\sum_j av_j \o w_j$ for any $a\in U$.

For any $n\geq N$, $W_n$ is included in $W_N W_{n-N}$, and one can write down 
$$w= \sum_j w_{j,N} w'_{j,n-N}$$
where $w_{j,N}\in W_N=R$ and $w'_{j,n-N} \in W_{n-N}$, so
$$d_l^N(a\o w)= \sum_j a\, \overline{w_{j,N}} \o w'_{j,n-N}.$$
Here $\overline{w_{j,N}}$ is the class of $w_{j,N}$ in $U^N$. But, in $U$, one has $\overline{w_{j,N}}= \varphi (w_{j,N})\in k$. Thus
$$d_l^N(a\o w)= \sum_j a\, \varphi(w_{j,N}) \o w'_{j,n-N}= a\o \sum_j \varphi(w_{j,N})w'_{j,n-N},$$so that 
$d_l^N = 1_U \o \varphi ^{1,N}$. Recall that $\varphi ^{1,N}:R\o V^{\o (n-N)}\rightarrow V^{\o (n-N)}$ denotes $\varphi \o 1_{V^{\o (n-N)}}$. Notice that $d_l$ is an $N$-differential if and only if $\varphi =0$.

Define analogously the $k$-$U$-linear map $d_r : W_n\o U \rightarrow W_{n-1}\o U$. Then $d_r^N= \varphi^{n-N+1,n}\o 1_U$, where $\varphi^{n-N+1,n}= 1_{V^{\o (n-N)}}\o \varphi$.

For the sake of simplicity, the $U$-$U$-linear maps $d_l\o 1_U$ and $1_U\o d_r$ will be also denoted by $d_l$ and $d_r$ respectively. Fix a primitive $N$-root of unity $q$ (we enlarge $k$ if necessary). Define 
$d:\ U\otimes W_{n}\otimes U \rightarrow U\otimes W_{n-1}\otimes U$
by $d=d_{l}-q^{n-1}d_{r}$. Explicitly we have
\begin{equation} \label{Ncomp}
\cdots  \stackrel{d_{l}-d_{r}}\longrightarrow U\otimes W_{N} \otimes U
\ \stackrel{d_{l}-q^{N-1}d_{r}}\longrightarrow \ \cdots  \stackrel{d_{l}-qd_{r}}\longrightarrow
U\otimes V \otimes U \stackrel{d_{l}-d_{r}}\longrightarrow U\otimes U.
\end{equation}  
\Blm
One has $d^N=0$, i.e., $(U\o W_n\o U)_{n\geq 0}$ endowed with $d$ is an $N$-complex.
\Elm
\Bdm
Since $d_l$ and $d_r$ are commuting, we have in $U\otimes W_{n}\otimes U$
$$\prod_{i=n-N}^{i=n-1} (d_{l}-q^{i}d_{r})=d_{l}^{N}-d_{r}^{N}=1_U\o(\varphi ^{1,N}-\varphi^{n-N+1,n})\o 1_U.$$
Proposition \ref{car4} shows that the PBW condition (3.5) reduces in our case to the following
$$(\varphi ^{1,N}-\varphi^{2,N+1})(W_{N+1})=0.$$
It implies that $(\varphi ^{1,N}-\varphi^{2,N+1})(W_{N+1}V)=(\varphi ^{2,N+1}-\varphi^{3,N+2})(VW_{N+1})=0$, hence $(\varphi ^{1,N}-\varphi^{3,N+2})(W_{N+2})=0$. A straightforward induction provides $(\varphi ^{1,N}-\varphi^{n-N+1,n})(W_n)=0$, so $d^N:\ U\otimes W_{n}\otimes U \rightarrow U\otimes W_{n-N}\otimes U$ is vanishing.
\qed
\\

The associated graded of the $N$-complex $(U\o W_n\o U)_{n\geq 0}$ is the Koszul $N$-complex of $A$ in $A$-grMod-$A$ as defined in \cite{bm:dual}. Therefore, since $A$ is Koszul and the filtration is exhaustive and bounded below, the adequate contraction $(U\o W_{\zeta(n)}\o U)_{n\geq 0}$ is exact in degrees $>0$. Recall that $\zeta : \mathbb{N} \rightarrow \mathbb{N}$, where
$$\zeta(2q)=qN, \ \zeta(2q+1)= qN+1,\ q\ \mathrm{integer} \geq 0.$$ 
Explicitly, the adequate contraction is the following complex
\begin{equation} \label{2comp}
\cdots  \stackrel{d^{N-1}}{\longrightarrow} U\otimes W_{N+1} \otimes U
\stackrel{d}{\longrightarrow} U\otimes W_{N} \otimes U \stackrel{d^{N-1}}{\longrightarrow}
U\otimes V \otimes U \stackrel{d}{\longrightarrow} U\otimes U.
\end{equation} 
Here $d=d_{l}-d_{r}$ and $d^{N-1}=d_{l}^{N-1}+d_{l}^{N-2}d_{r} + \cdots
+ d_{l}d_{r}^{N-2}+d_{r}^{N-1}$ make sense on the genuine $k$. Together with $\mu : U\o U \rightarrow U$, $(U\o W_{\zeta(n)}\o U)_{n\geq 0}$ is a resolution of $U$ in the category of filtered bimodules $U$-filtMod-$U$, which is called the \emph{Koszul resolution} of $U$.

If $n$ is odd, one has
\begin{align*} 
&d(a\o\sum_i v_1^i\o \cdots \o v_{\zeta(n)}^i\o b)\\
&=\sum_i \left(av_1^i\o v_2^i\o \cdots \o v_{\zeta(n)}^i\o b- a\o v_1^i\o \cdots \o v_{\zeta(n-1)}^i\o v_{\zeta(n)}^i b \right),
\end{align*}
and if $n$ is even:
\begin{align*} 
&d^{N-1}(a\o\sum_i v_1^i\o \cdots \o v_{\zeta(n)}^i\o b)\\
&=\sum_i \big(av_1^i \cdots v_{N-1}^i \o v_N^i \o \cdots \o v_{\zeta(n)}^i \o b + av_1^i \cdots v_{N-2}^i \o v_{N-1}^i \o \cdots \o v_{\zeta(n-1)}^i \o v_{\zeta(n)}^i b\\
&+ \cdots + a\o v_1^i \o \cdots \o v_{\zeta(n-N+1)}^i \o v_{\zeta(n-N+2)}^i \cdots v_{\zeta(n)}^i b \big).
\end{align*}

Denoting the Koszul resolution of $U$ by $\mathcal{K}(U)$, $\mathcal{K}(U)$ is a $k$-split ($U^e\otimes$)-projective resolution if $k$ is semi-simple \cite{weib:homo}. For any $U$-$U$-bimodule $M$, let us define the Hochschild homology and cohomology of the $k$-$k$-algebra $U$ (here $k$ may be non-commutative!) by 
$$HH_{\bullet}(U,M)= \mathrm{Tor}_{\bullet}^{U^e/k}(M,U), \ \ HH^{\bullet}(U,M)= \mathrm{Ext}^{\bullet}_{U^e/k}(U,M).$$
So if $k$ is semi-simple, one has
$$HH_{\bullet}(U,M)=H_{\bullet}(M\o_{U^e}\mathcal{K}(U)), \ \ HH^{\bullet}(U,M)=H^{\bullet}(\mathrm{Hom}_{U^e}(\mathcal{K}(U),M)).$$

Let us specialize to the situation of Theorem \ref{maint}, i.e. $N=p$ and $U=H_{\psi}$ which is supposed Koszul. Here $W_n$ is identified to $\wedge^{n}_K E$. So for any $v_1, \ldots v_n$ in $V$, Alt$(v_1, \ldots , v_n)$ is identified to the wedge product $v_1\wedge\ldots \wedge v_n$ \emph{defined over} $K=k[\Gamma ]$. In the following formulas, the tensor products are also \emph{defined over} $K$. So, if $n$ is odd:
\begin{align*}
&d(a\o v_1\wedge  \cdots \wedge v_{\zeta(n)}\o b)\\ 
&=\sum_{j=1}^{\zeta(n)}\big((-1)^{j-1} av_j\o v_1 \wedge \cdots \wedge \hat{v_j} \wedge \cdots \wedge v_{\zeta(n)}\o b\\
&-(-1)^{\zeta(n)-j}  a\o v_1 \wedge \cdots \wedge \hat{v_j} \wedge \cdots \wedge v_{\zeta(n)}\o v_j b \big).
\end{align*}
Note that when $p$ is even, $(-1)^{\zeta(n)-j}=(-1)^{j-1}$. Now if $n$ is even:
\begin{align*} 
&d^{p-1}(a\o v_1\wedge  \cdots \wedge v_{\zeta(n)}\o b)=\sum_{1\leq j_1<\cdots <j_{p-1}\leq \zeta(n)}\\ \\
&\big((-1)^{j_1-1+\cdots +j_{p-1}-(p-1)} av_{j_1}\ldots v_{j_{p-1}}\o v_1\wedge  \cdots  \hat{v}_{j_1}  \cdots  \hat{v}_{j_{p-1}}  \cdots \wedge  v_{\zeta(n)} \o b\\ \\
& +(-1)^{j_1-1+\cdots +j_{p-2}-(p-2)+ \zeta(n)-j_{p-1}} av_{j_1}\ldots v_{j_{p-2}}\o v_1\wedge  \cdots  \hat{v}_{j_1}  \cdots  \hat{v}_{j_{p-1}}  \cdots \wedge  v_{\zeta(n)} \o v_{j_{p-1}}b\\ \\
&+ \cdots + (-1)^{\zeta(n)-j_{p-1}+ \cdots + \zeta(n)-j_{1}-(p-2)} a\o v_1 \wedge  \cdots  \hat{v}_{j_1}  \cdots  \hat{v}_{j_{p-1}}  \cdots \wedge  v_{\zeta(n)} \o v_{j_1}\ldots v_{j_{p-1}}b \big),
\end{align*}
which is reduced to the following when $p$ is even:
\begin{align*} 
&d^{p-1}(a\o v_1\wedge  \cdots \wedge v_{\zeta(n)}\o b)=\sum_{1\leq j_1<\cdots <j_{p-1}\leq \zeta(n)}(-1)^{j_1+\cdots +j_{p-1}+\frac{p}{2}}\\ \\
&\big(av_{j_1}\ldots v_{j_{p-1}}\o v_1\wedge  \cdots  \hat{v}_{j_1}  \cdots  \hat{v}_{j_{p-1}}  \cdots \wedge  v_{\zeta(n)} \o b\\ \\
& - av_{j_1}\ldots v_{j_{p-2}}\o v_1\wedge  \cdots  \hat{v}_{j_1}  \cdots  \hat{v}_{j_{p-1}}  \cdots \wedge  v_{\zeta(n)} \o v_{j_{p-1}}b\\ \\
&+ \cdots - a\o v_1 \wedge  \cdots  \hat{v}_{j_1}  \cdots  \hat{v}_{j_{p-1}}  \cdots \wedge  v_{\zeta(n)} \o v_{j_1}\ldots v_{j_{p-1}}b \big),
\end{align*}
If $p=2$, the formulas for $n$ odd and $n$ even are the same and we recover in this case the differential as defined in \cite{eg:reflec}, formula (2.6).


\begin{thebibliography}{99}

\bibitem{as:regular} M. Artin, W. F. Schelter, Graded algebras of global dimension 3, \emph{Adv.
Math.} \textbf{66} (1987), 171-216.
\bibitem{atv:grot} M. Artin, J. Tate, M. Van den Bergh, Some algebras associated to automorphisms
of elliptic curves, The Grothendieck Festschrift, vol 1,
Birkh\"{a}user, Basel, 1990.
\bibitem{atv:mod} M. Artin, J. Tate, M. Van den Bergh, Modules over
  regular algebras of dimension 3, \emph{Invent.
math.} \textbf{106} (1991), 335-388.
\bibitem{bgso:kdp} A. A. Beilinson, V. Ginzburg, W. Soergel, Koszul duality patterns in
representation theory, \emph{J. Am. Math. Soc.} \textbf{9} (1996), 473-527.
\bibitem{br:dup} G. Benkart, T. Roby, Down-up algebras, \emph{J. Algebra}
\textbf{209} (1998), 305-344.
\bibitem{rb:nonquad} R. Berger, Koszulity for nonquadratic algebras, \emph{J. Algebra}
\textbf{239} (2001), 705-734.
\bibitem{rb:nonquad2} R. Berger, Koszulity for nonquadratic algebras II, math.QA/0301172.
\bibitem{bdvw:homog} R. Berger, M. Dubois-Violette, M. Wambst, Homogeneous algebras,
\emph{J. Algebra} \textbf{261} (2003), 172-185.  
\bibitem{bm:dual} R. Berger, N. Marconnet, Koszul and Gorenstein properties for homogeneous
algebras, \emph{Alg. Rep. Theory}, to appear.
\bibitem{bou:algco} N. Bourbaki, \emph{Alg\`ebre commutative}, Chap. 1-4, \'El\'ements de
math\'ematique, Masson, 1985.
\bibitem{bg:pbw} A. Braverman, D. Gaitsgory, Poincar\'e-Birkhoff-Witt
  theorem for quadratic algebras of Koszul type, \emph{J. Algebra}
\textbf{181} (1996), 315-328.
\bibitem{eg:reflec} P. Etingof, V. Ginzburg, Symplectic reflection
  algebras, Calogero-Moser space, and deformed Harish-Chandra homomorphism, \emph{Invent. math.} \textbf{147} (2002), 243-348.        
\bibitem{good:von} K. R. Goodearl, \emph{Von Neumann regular rings}, Pitman, 1979.
\bibitem{pp:quad} A. Polishchuk, L. Positselski, Quadratic algebras,
  preprint, 1994.
\bibitem{weib:homo} C. A. Weibel, \emph{An introduction to homological algebra}, Cambridge University
Press, 1994.

\end{thebibliography}
\end{document}